\documentclass[10pt]{article}

\usepackage{latexsym}
\usepackage{amsthm, amsmath, amssymb,latexsym,xypic, color}

\newtheorem{theorem}{Theorem}[section]

\newtheorem{lemma}[theorem]{Lemma}

\newtheorem{proposition}[theorem]{Proposition}
\newtheorem{remark}[theorem]{Remark}

\newcommand{\nc}{\newcommand}
\nc{\calH}{{\mathcal H}}
\nc{\calS}{{\mathcal S}}
\newcommand{\im}{\mbox{ Im }}
\nc{\uR}{\underline{\mathbb R}}
\nc{\uC}{\underline{\mathbb C}}
\nc{\bH}{{\mathbb H}}
\nc{\bA}{{\mathbb A}}
\nc{\bG}{{\mathbb G}}
\nc{\bC}{{\mathbb C}}
\nc{\bO}{{\mathbb O}}
\nc{\bI}{{\mathbb I}}
\nc{\bB}{{\mathbb B}}
\nc{\bY}{{\mathbb Y}}
\nc{\bK}{{\mathbb K}} 
\nc{\bX}{{\mathbb X}}
\nc{\bS}{{\mathbb S}}
\nc{\bE}{{\mathbb E}}
\nc{\bF}{{\mathbb F}}
\nc{\bZ}{{\mathbb Z}}
\nc{\bQ}{{\mathbb Q}}
\nc{\bN}{{\mathbb N}}
\nc{\bP}{{\mathbb P}}
\nc{\bL}{{\mathbb L}}
\nc{\bM}{{\mathbb M}}
\nc{\bT}{{\mathbb T}}
\nc{\bW}{{\mathbb W}}
\nc{\bU}{{\mathbb U}}
\nc{\bD}{{\mathbb D}}
\nc{\bJ}{{\mathbb J}}
\nc{\bV}{{\mathbb V}}
\nc{\bbZ}{{\mathbb Z}}
\nc{\bR}{{\mathbb R}}
\nc{\co}{{\nabla}}
\nc{\cu}{{\overline{\nabla}}}
\nc{\fr}{{\rightarrow}}

\begin{document}
                                %
                                %
                                %
\title{A note on spaces of symmetric matrices} %
\author{Andrea Causin and Gian Pietro Pirola  \footnote{Partially supported by
1) PRIN 2005 {\em ``Spazi di moduli e
teorie di Lie"}; 
2) Indam (GNSAGA);
3) Far 2006 (PV):{\em ``Variet\`{a} algebriche, calcolo
algebrico, grafi orientati e topologici"}.}}
\date{}     

                                %
                                %
\maketitle
                                %
                                %

                                %
                                %


                                %
                                %
                                %
                                %
                                %
 \begin{abstract} {\em 
\noindent We calculate the maximal dimension of linear spaces of symmetric and hermitian matrices with given high rank generalizing a well-known result of Adams et al.} \vskip3mm
\vskip 1mm
 {\setlength{\baselineskip}{0.8\baselineskip}
 \noindent {\scriptsize {\bf AMS (MOS) Subject 
 Classification:} {\em 15A30 (55N15).} }\\
  \noindent {\scriptsize {\bf Key words:} {\em  Symmetric matrices, K-theory, homotopy of classical groups.}} \par}

\end{abstract}

 \section*{Introduction}
Let $X$ denote a set of matrices over a field; we say that $V$ is a $k-$space in $X$ whenever $V\subset X$ is a real vector space whose nonzero elements have rank $k$. A natural problem in this context is to determine the maximal dimension $d_X(k) $ of a $k-$space in a given $X$.
For real invertible matrices, the  answer  has been given by Adams  by determining the maximal number of independent vector fields on a sphere \cite{adams1}. His work provides the keystone for studying  interesting cases. In  \cite{adams2} the solution is given for invertible symmetric real matrices, invertible complex and quaternionic matrices, as well as for their hermitian relative cases.
Many subsequent researches on the above problem and  its generalizations (e.g. to the case of matrices with bounded rank) has been done; see, among others  \cite{yeung}, \cite{rees1}, Friedland {\em et al.} \cite{pacific,oddeg,crossing}. 
In particular, in \cite{rees1}  $X$ is the set of the real matrices with fixed rank 
 and in \cite{pacific,crossing} the setup in terms of nonlinear problems over spheres is explicitly realized.

In the present paper we study the real symmetric and hermitian matrices. These appear  in
several different areas, e.g. hyperbolic system of differential equation, spectral problems and
cohomology of K\"ahler varieties \cite {pacific,adams2,causin}.   

For real $r$, define the Radon-Hurwitz numbers $\rho(r) = 2^c+8d$ and $\rho_\bC(r)=2(c+4d)+2$ when $r=2^{c+4d}(2a+1)$, with $a,c$ and $d$ integers, $0\leq c\leq 3$; $\rho_\bC(r)=\rho(r)=0$ otherwise. 
Moreover set  $\sigma(n,h)=\max\{\rho(\frac{h}{2}+j) \mbox{ with } 0\leq j \leq n-h \}$ and $\sigma_\bC$ in a similar fashion using $\rho_\bC$. We prove:

\medskip
\noindent {\bf Theorem 1}\\
{\em Let  $X$ be the set of $n\times n$ real symmetric matrices and $0\leq s\leq 2;$  then 
\begin{equation}\label{bounds} \sigma(n,n-s) \leq d_X(n-s)\leq\sigma(n,n-s)+1. \end{equation}  When $\sigma(n,n-s)=\rho(\frac{n-s}{2})$, the upper bound is attained. \\ 
If $s=1$, the lower bound is optimal when $\frac{n+1}{2}=2,$ $2^{2+4d}\gamma$ or $ 2^{3+4d}\gamma$, where $d$ is an arbitrary integer and $\gamma$ is an odd integer. In this case, $\sigma(n,n-1)=\rho(\frac{n+1}{2})$.}
\medskip

\noindent {\bf Theorem 2}\\
{\em Let $X$ be the set of $n\times n$ complex hermitian matrices and $s=0,1$; then 
\begin{equation}\label{boundsh} \sigma_\bC(n,n-s) \leq d_X(n-s)\leq\sigma_\bC(n,n-s)+1. \end{equation}  If $\sigma_\bC(n,n-s)=\rho_\bC(\frac{n-s}{2})$, the upper bound is attained, otherwise the lower bound is optimal.}
\smallskip

In the case $s=0,$ our theorems provide a new proof of some of the results in \cite{adams2},  and in the case $s=1$ Theorem 1 improves the estimate given in \cite{pacific}. 

We notice  that Theorems $1$ and $2$ can be rephrased in terms of maps from spheres to spaces of matrices (see e.g. \cite{crossing}); that is there exists an odd continuous map $\phi:S^d\rightarrow X$ if and only if $d< d_X.$  

\smallskip 

The scheme of the paper is the following: in the first section we prove Theorem $1;$ the proof is divided in three parts showing respectively the upper and the lower bounds, and their optimality in the stated cases; the second section is devoted to prove Theorem 2.

\smallskip

It is a pleasure to thank Prof.~S.~Friedland for the helpful papers he brought to our attention; in particular, a keypoint of our work relies on the ideas of \cite{crossing}. We are grateful to Margherita for her suggestions and for the help she provided us.
 
\section{Symmetric matrices}
\subsection*{Upper bound}

Let  $X=\calS_n^k,$ where $\calS_n^k$ is the set of $n-$square real symmetric matrices of rank $k.$ 
In this paragraph we will give a prove of the upper bound: 
\begin{proposition}\label{upper} 
The following inequality holds: 
$$d_X(n-s)\leq\sigma(n,n-s)+1.$$ 
\end{proposition}
We notice that  the Proposition \ref{upper} is equivalent to the inequality $$d\leq\sigma(n,n-s) \quad\mbox{\em  for a sphere $S^d\subset \calS_n^{n-s}$  or  a projective space $\bP^d\subset \bP(\calS_n^{n-s})$}.$$ This is the form we will be referring to.\\
We also note that $S^d\subset \calS_n^k$ implies $k$ even or $d=0;$ indeed, if $d>0$ there is a path (of constant rank matrices) in $S^d$ connecting any matrix $A$ to $-A$ and this forces the signature of $A$ to be $(\frac{k}{2},\frac{k}{2}).$

As in \cite{rees1}, over $\bP^d\subset \bP(\calS_n^{k})$ with $d>0$, we can construct the exact sequence of bundles:

\begin{equation}\label{evaluationg}
0\longrightarrow K \longrightarrow  \uR^n \stackrel{E}{\longrightarrow} H^n\longrightarrow C\longrightarrow 0
\end{equation} where $H$ is the hyperplane nontrivial line bundle.
The central map is given by $E([A],v)=([A], Av)$ and, since all matrices $A$ have constant rank $k$, its kernel defines a kernel bundle $K$ and a cokernel bundle $C$;  moreover, the isomorphism  $K\oplus H^n \simeq \uR^n \oplus C$ holds.

Denote by $\pi:S^d\rightarrow \bP^d$ the quotient of the multiplication by $-1$ and consider the pullback of sequence (\ref{evaluationg}) via $\pi.$

Now, we need to show two preliminary lemmas. 
\begin{lemma}\label{epiu} There exist isomorphic bundles $E^+$, $E^-$ over $S^d$ such that $\pi^\ast K\oplus E^+\oplus E^- = \uR^n$. Their rank is $\frac{k}{2}$.
\end{lemma}
\begin{proof} Let us say that an eigenvector is positive (resp.~negative) if it is relative to a positive (resp.~negative) eigenvalue. Let $E^+$ (resp. $E^-$) be the bundle whose fiber over a matrix is the span of its positive (resp.~negative) eigenvectors. Clearly, if $v$ is a positive eigenvector for $A$, it is negative for $-A$, hence multiplication by $-1$ on $S^d$ lifts to an automorphism of $E^+\oplus E^-$ interchanging the summands.
\end{proof}

\begin{lemma}\label{sommadiff} If $E^+$ is trivial, then $E^+\oplus E^-$ is isomorphic to $\pi^\ast (\uR\oplus H)^{\frac{k}{2} }$.
\end{lemma}
\begin{proof} Since $E^+$ is trivial, we can choose a basis (at any point) $v_1^+,\dots ,v_{\frac{k}{2}}^+$ of global sections of $E^+,$ and we construct the corresponding sections of $E^-,$ $v_1^-,\dots ,v_{\frac{k}{2}}^-,$  by setting $v_i^-(A)=-v_i^+(-A).$ Therefore, we define new sections for $E^+\oplus E^-$ as follows: $$r_i = v_i^+ - v_i^- \quad \mbox{ and } \quad h_i = v_i^+ + v_i^-.$$ These new sections decompose $E^+ \oplus E^-$ as a sum of $k$ trivial line bundles $L_j$. Since the sections $r_i$ are invariant under the action of $-1$ in $S^d$, their corresponding line bundles are the pullback of  $\uR\rightarrow\bP^d;$ on the other hand, the $h_i$ are anti-invariant with respect to the same action, and this shows that the remaining line bundles are the pullback of $H\rightarrow \bP^d$.
\end{proof}

To complete the proof of the Proposition \ref{upper}, it is now sufficient to notice that any bundle over $S^d$ becomes trivial when it is restricted to $S^{d-1}$ (equator), 
since $S^d$ minus a point is contractible. 
Applying the above Lemma \ref{sommadiff} to this restriction,  we get the relation 
\begin{equation}\label{relazione} K\oplus \uR^{\frac{k}{2}}\oplus H^{\frac{k}{2}}=\uR^n \quad \mbox{over $\bP^{d-1}$.}  \end{equation}  
Then, the estimate of the Proposition \ref{upper} is a consequence of the following well-known fact (cf. \cite{adams1,rees1}): 

\medskip
{\em The reduced ring of real K-theory $\tilde K(\bP^{d-1})$ is the ring $\bZ[\mu]$ with relations $\mu^2=-2\mu $ and $2^f\mu=0$, where $\mu = [H]-[\uR]$ and $f$ is the number of integers $s\equiv 0,1,2$ or $4$ $\mod8$ such that $0< s < d$. In particular, for any integer $m,$ $m\mu=0$ implies $d\leq\rho(m)$.}

\smallskip 

From this, we immediately get:
\begin{enumerate}
\item[$\ast$] $s=0:$ we get $K=\underline0$ and $k=n$; the relation (\ref{relazione}) entails $\frac{n}{2}\mu=0$ i.e. $d\leq\rho(\frac{n}{2})=\sigma(n,n);$
\item[$\ast$] $s=1$ (that is $k=n-1$): $K$ is either $\uR$ or $H$; in the first case we get $\frac{n-1}{2}\mu=0;$ in the second $\frac{n+1}{2}\mu=0$; this implies $d\leq\sigma(n,n-1)$; 
\item[$\ast$] $s=2:$ $K$ is $\uR^2,$  $\uR\oplus H$ or $H^2$ (cf. e.g. \cite{adams3}), thus respectively $\frac{n-2}{2}\mu,\ \frac{n}{2}\mu$ or $\frac{n+2}{2}\mu$ is zero, and $d\leq\sigma(n,n-2).$
\end{enumerate}

\subsection*{Lower bound} 
In this paragraph we prove the lower bound in Theorem 1 by showing that 
there are suitable spaces of matrices of the stated dimension.
The basic brick for constructing all examples we need is the fact (see \cite{adams2}) that: {\em there exists a $\rho(m)-$dimensional space $V_m$ of invertible $m\times m$ matrices such that any nonzero $A\in V_m$ verifies $^t\!AA=y^2 I$, for some real $y\neq 0$.}

Assume firstly $s=0$. It is possible to find a  space $W_n^0$ of dimension $\rho(\frac{n}{2})+1$ of $n\times n$ symmetric and invertible (excepted $0$) matrices;  it is the example provided in~\cite{adams2}:
\begin{equation}\nonumber 
\begin{pmatrix}  xI & A \\ ^t\!A & - xI \end{pmatrix} \quad A \in V_{\frac{n}{2}}, \, x\in \bR.
\end{equation}

If $s=1$ we construct two spaces of $n\times n$ matrices. The first one is straightforward:
\begin{equation}\nonumber 
\begin{pmatrix} 0 & 0 \\ 0 & M \end{pmatrix}, \quad M\in W_{n-1}^0;
\end{equation}
clearly, it has dimension $\rho(\frac{n-1}{2})+1.$ The second one is the space of the matrices \begin{equation}\nonumber 
\begin{pmatrix} 0 & B \\ ^t\!B & 0 \end{pmatrix}
\end{equation}
where $B$ is obtained from $A\in V_{\frac{n+1}{2}}$ by deleting the last row. Such $B$ are of maximal rank and form a space of dimension equal to $\dim V_{\frac{n+1}{2}}=\rho(\frac{n+1}{2});$ indeed, if this does not hold, there would exist $A\in V_{\frac{n+1}{2}}$ with all rows $0$ except the last one, which gives a contradiction. Constructed these two spaces, for each $n$ we denote by $W_{n}^1$ 
the one with maximal dimension: then  
$$ \dim W_{n}^1=  \max\Big\{\rho\big(\frac{n-1}{2}\big)+1,\, \rho\big(\frac{n+1}{2}\big)\Big\}.$$ 

Finally, assume $s=2$. As in the previous case, we construct either the space of matrices $\begin{pmatrix} 0 & 0 \\ 0 & M \end{pmatrix}$ with $M\in W_{n-1}^1$, whose dimension is $\max\{\rho(\frac{n-2}{2})+1,\, \rho(\frac{n}{2})\}$, or the space of $\begin{pmatrix} 0& B\\ ^t\!B & 0\end{pmatrix}$ with $B$ obtained from $A\in V_{\frac{n+2}{2}}$ by deleting the two last rows.

\begin{remark} The spaces constructed above show that when $\sigma(n,n-s)=\rho(\frac{n-s}{2}),$ then the upper bound in (\ref{bounds}) is reached, concluding the proof of the lower bound. 
\end{remark}

\subsection*{Optimality  of lower bound}
Here we complete the proof of the Theorem $1,$ by showing the last statement. 
The case $n=3$ is shown in \cite{pacific,crossing}; it remains to prove the following:  
\begin{proposition} \label{optimality} 
Assume that $\rho(\frac{n+1}{2})=  4+8d$ or $ 8+8d$ for some integer $d;$ then, $d_X(n-1)=\sigma(n,n-1) = \rho(\frac{n+1}{2}).$
\end{proposition}

Thanks to the lower bound inequality, it is sufficient to prove that the case $d_X(n-1)=\sigma(n,n-1)+1=\rho(\frac{n+1}{2})+1$ does not hold.  
We show this by contradiction; assume $d_X(n-1)=\sigma(n,n-1)+1=\rho(\frac{n+1}{2})+1$ and denote this number by $r+1$. The Lemma \ref{epiu} gives isomorphic bundles $E^{\pm}$ of rank $\frac{n-1}{2}$ over $S^r$. If we can show that these bundles are trivial, then the  Lemma \ref{sommadiff} would imply the relation $$K\oplus(\uR\oplus H)^\frac{n-1}{2}=\uR^n \mbox{ over } \bP^r.$$ 
This gives the contradiction $r+1\leq\sigma(n,n-1)=r.$ 

It remains to prove that, under the hypothesis $d_X(n-1)=\sigma(n,n-1)+1=\rho(\frac{n+1}{2})+1$, the bundles $E^{\pm}$ are trivial. 

Recall (see \cite{husemoller} for what follows) that real bundles of rank $k$ over $S^r$, with $r\geq2$, are classified up to isomorphism by the homotopy groups $\pi_{r-1}(SO(k))$; there is a natural inclusion $SO(k)\subset SO(k+1)$  whose induced map $i$ on the $m-$th homotopy groups is an isomorphism if $k>m+1$; in this case, Bott periodicity holds: $\pi_{m}(SO(k))=\pi_{m+8}(SO(k))$; moreover, if $A$ and $B$ are maps representing bundles $F $ and $G$, the map representing $F\oplus G$ is $\begin{pmatrix}A & 0 \\ 0 & B \end{pmatrix} = \begin{pmatrix}A & 0 \\ 0 & I \end{pmatrix} \begin{pmatrix} I & 0 \\ 0 & B \end{pmatrix}= i^{rkG}A + i^{rkF}B \in\pi_{r-1}(SO(rkF+rkG))$. 

\smallskip

Now, we can show the following two lemmas, corresponding to the cases we are dealing with.  
\begin{lemma}\label{minimo} If $r=4+8d$, then $E^\pm$ are trivial bundles.
\end{lemma}
\begin{proof} We show that the map $i^{\frac{n-1}{2}}:\pi_{r-1}(SO(\frac{n-1}{2}))\rightarrow \pi_{r-1}(SO(n-1))$ is injective and the target group is isomorphic to $\bZ$. This will conclude the proof, since  $E^+\oplus E^-=E^+\oplus E^+ =\uR^{n-1}$ and if $e$ represents $E^+$ we will get $2i^{\frac{n-1}{2}} e=0$.

Observe that $r=4+8d$ is equivalent to $\frac{n+1}{2}=2^{2+4d}\gamma$ with $\gamma$ odd, hence the above map is $$ \pi_{3+8d}(SO(16^d4\gamma-1))\longrightarrow \pi_{3+8d}(SO(n-1))$$ and is a composition of isomorphisms provided $16^d4\gamma-1>3+8d+1$ that is $\gamma\neq1$ and $d\neq0$; moreover, all those groups are isomorphic to $\bZ$ thanks to Bott periodicity and the fact that $\pi_3(SO(k))=\bZ$ stably.

Then, take $d=0$ and $\gamma=1$.  The corresponding map is the composition $$ \pi_3(SO(3))\stackrel{i}{\rightarrow}\pi_3(SO(4))\stackrel{j}{\rightarrow}\pi_3(SO(5)) \rightarrow \pi_3(SO(6));$$ the last arrow is a stable isomorphism $\bZ\rightarrow \bZ$, thus we only need to show that $ji$ is not zero. Computing the exact homotopy sequence of $SO(3)\rightarrow SO(4)\stackrel{p}{\rightarrow} S^3$ shows that $\pi_3(SO(3))=\bZ$, $\pi_3(SO(4))=\bZ\oplus\bZ$ and $i$ is injective. Moreover,  $\im i=\ker p_\ast$ and $\ker j=\im \partial$ where $\partial$ is the injective boundary in the sequence $\pi_4(S^4)\stackrel{\partial}{\rightarrow}\pi_3(SO(4))\stackrel{j}{\rightarrow}\pi_3(SO(5))\rightarrow\pi_3(S^4)=0$.
 In (\cite[thm.~10.4]{husemoller}) it is shown that $\im \partial$ is generated by the characteristic map $c:S^3\rightarrow SO(4)$ of the principal bundle  associated to the tangent bundle of $S^4$. 
 It is also shown (\cite[thm.~10.1]{husemoller}) that the composition $pc:S^3\rightarrow S^3$ has degree $2$. This forces  $\pi_4(S^4)\stackrel{\partial}{\rightarrow}\pi_3(SO(4))\stackrel{p_\ast}{ \rightarrow}\pi_3(S^3)$ to be the multiplication by $2$, indeed $p_\ast\partial([id])=p_\ast([c])=[pc]$. Hence $\im \partial \cap \ker p_\ast =\{0\}$, so $ji$ is not zero.
\end{proof}

\begin{lemma} If $r=8+8d$, then $E^\pm$ are trivial bundles.
\end{lemma}

\begin{proof} We  argue as in the previous lemma. Now we deal with maps $$\pi_{7+8d}(SO(16^d8\gamma-1))\longrightarrow \pi_{7+8d}(SO(n-1)) $$ that fall in the range of  stable inclusion of homotopy groups when $d\neq0$ and $\gamma\neq 1$, hence  they all are isomorphisms $\bZ\rightarrow\bZ$. The only case left is $\pi_7(SO(7))\rightarrow\pi_7(SO(14))$ which reduces to determine $\pi_7(SO(7))\rightarrow\pi_7(SO(9)),$ but this is done exactly as before (cf. \cite{hatcher,husemoller}). 
\end{proof}

\begin{remark} 
When $n\equiv 3\mod 4$, we get $\sigma(n,n-1)=\rho(\frac{n+1}{2})$ and this number is written $2^c+8d$ with $0\leq c\leq 3$. Proposition \ref{optimality} says that if $c=2$ or $3$, then the lower bound of Theorem 1 is the exact estimate of $d_X(n-1)$.  
On the other hand, the methods developed in this section can not be used to decide the optimality of such bound when $c=0$ or $1$; indeed, in those cases, the stable homotopy groups $\pi_{r-1}(SO(k))$ are cyclic of order 2. 
Also note that a statement similar to Proposition \ref{optimality} could be proved for $s=2$.
\end{remark}

\section{Hermitian matrices}
The outline of the proof of Theorem 2 is essentially the same of Theorem 1. 
We only remark the adapted steps.

\smallskip 

\noindent {\em Upper bound.} The calculations done in the previous section can be adapted simply
 using complex bundles instead of real ones, since hermitian matrices have real eigenvalues and 
 there is only one nontrivial complex line bundle $H_\bC$ over the real projective space. 
 Moreover, the ring of complex K-theory $K_\bC(\bP^{d-1})$   is generated by $\nu=[H_\bC]-[\uC]$ and provides the implication $m\nu=0\Rightarrow d\leq\rho_C(m)$. 

\smallskip 

\noindent {\em Lower bound.} In \cite{adams2}, it is shown that there exist  $\rho_\bC(m)$ complex  $m\times m$ matrices  ``whose real linear combinations are nonsingular"; then, we can follow exactly the construction we did in the real case.

\smallskip 

\noindent {\em Optimality.} Clearly, the only case we have to consider is  $s=1$: 
if $n$ is even, there is nothing to prove. 
If $n\equiv1\mod 4$, then $\sigma_\bC(n,n-1)=\rho_\bC(\frac{n-1}{2})$ and the bound is reached by explicit examples.
When $n\equiv3\mod 4$, then $\sigma_\bC(n,n-1)=\rho_\bC(\frac{n+1}{2});$ we denote this number by $r,$ and we show that the
 upper bound can never be attained. Suppose by contradiction that it is attained; then, with
 the same argument of Proposition \ref{optimality}, we need to prove that the bundles $E^\pm$ are trivial.
We have to study the homotopy maps  $\pi_{r-1}(SU(k))\rightarrow \pi_{r-1}(SU(2k)).$  These maps are  isomorphisms for $n\neq 3,$ as can be seen by computing the homotopy sequences of $SU(m)\rightarrow SU(m+1)\rightarrow S^{2m-1}.$   Since $r$ is always even, complex Bott periodicity ensures that these groups are isomorphic to $\bZ.$ Finally, if $n=3$, line bundles on the $4-$sphere are trivial since $\pi_3(SU(1))=0$.

\noindent
 {\sc Andrea Causin}\\
Dipartimento di Matematica, Universit\`a {\em``La Sapienza"} di Roma\\
P.le Aldo Moro 2, 00185 Roma, Italia\\
{\tt causin@mat.uniroma1.it}\\

\noindent {\sc Gian Pietro Pirola}\\
Dipartimento di Matematica, Universit\`a di Pavia\\
via Ferrata 1, 27100 Pavia, Italia\\
{\tt gianpietro.pirola@unipv.it}

\end{document}